\documentclass{amsart}

\usepackage{amsmath,amsthm}
\usepackage{amsfonts,amssymb}

\newtheorem{thm}{Theorem}[section]
\newtheorem{cor}[thm]{Corollary}
\newtheorem{lem}[thm]{Lemma}
\newtheorem{prop}[thm]{Proposition}

\theoremstyle{remark}
\newtheorem{rem}{Remark}[section]

 \def\CL{{\mathcal L}}
\newcommand{\cI}{{\mathrm{I}}}
\newcommand{\cII}{{\mathrm{I\!I}}}

\def\a{{\alpha}} 
 \def\b{{\beta}}
 
 \def\k{{\kappa}}
 
 \def\l{{\lambda}}
 
 \def\o{{\omega}}
 \def\s{{\sigma}}
 \def\la{{\langle}}
 \def\ra{{\rangle}} 
 
\def\na{{\nabla}}
 
 \def\CH{{\mathcal H}}
 \def\CJ{{\mathcal J}}
 \def\CD{{\mathcal D}}
 \def\CV{{\mathcal V}}
 \def\RR{{\mathbb R}}

 \def\proj{\operatorname{proj}}

\newcommand{\wh}{\widehat}

\begin{document}

\title
{Sobolev Orthogonal Polynomials Defined via Gradient on the Unit Ball }

\author{Yuan Xu}
\address{Department of Mathematics\\ University of Oregon\\
    Eugene, Oregon 97403-1222.}\email{yuan@math.uoregon.edu}

\date{\today}
\thanks{The work was supported in part by NSF Grant DMS-0604056}

\begin{abstract}
An explicit family of polynomials on the unit ball $B^d$ of $\RR^d$ is 
constructed, so that it is an orthonormal family with respect to the inner 
product 
$$
 \langle f,g \rangle
    = \rho \int_{B^d}\nabla f(x)\cdot \nabla g(x) dx + \CL (fg), 
$$
where $\rho >0$, $\nabla$ is the gradient, and $\CL(fg)$ is either the 
inner product on the sphere $S^{d-1}$ or $f(0)g(0)$. 
\end{abstract}

\maketitle

\section{Introduction}\label{Introduction}
\setcounter{equation}{0}

In a problem related to dwell time for polishing tools in fabricating optical 
surfaces, it is important to control the gradient as it is related to mechanical 
accelerations of the polishing machines  (\cite{GF}). One then asks the question if 
there is a Parseval type relation for the square of the gradient. This problem
leads us to the study of the Sobolev type orthogonal polynomials with respect
to the inner product defined in terms of gradient on the unit ball,
\begin{equation}\label{ip1}
  \langle f,g \rangle_\cI := \rho \int_{B^d} 
     \na  f(x) \na  g(x) dx +  \int_{S^{d-1}} f(x)g(x) d\o(x),
\end{equation}
where $\rho > 0$,  $B^d$ and $S^{d-1}$ are the unit ball and the unit sphere
of $\RR^d$, respectively. Throughout this paper we write
$$
    \na f  \na g = \na f \cdot \na g : = \sum_{i=1}^d 
         \frac{\partial f}{\partial_i}  \frac{\partial g}{\partial_i} .
$$
Our main concern is the gradient term, the second 
term is added to make the inner product well defined. Clearly we can add
some other terms to the gradient part to make the inner product well defined.
We consider another choice
\begin{equation}\label{ip2}
  \langle f,g \rangle_{\cII} := \rho \int_{B^d} 
     \na  f(x) \na  g(x) dx +  f(0)g(0),
\end{equation}
where again $\rho>0$. Our goal is to find a complete system of orthonormal 
polynomials with respect to these inner products and explore their properties. 

The orthogonal polynomials with respect to these inner products are called 
the Sobolev orthogonal polynomials. Although there are many papers in 
the literature dealing with Sobolev orthogonal polynomials of one variable
(see, for example, \cite{G} and the references therein), as far as we know
\cite{X06} is the only paper that discusses Sobolev orthogonal polynomials 
of several variables. The study in \cite{X06} starts from a question on the
numerical solution of the nonlinear Poisson equation $-\Delta u = f(\cdot, u)$ 
on the unit disk with zero boundary conditions (see \cite{A}), which asks for
an explicit orthogonal basis for the inner product defined by
\begin{equation}\label{innerDelta}
  \langle f,g \rangle_\Delta := \frac{1}{\pi}\int_{B^d} 
        \Delta[(1-\|x\|^2) f(x)]\Delta[(1-\|x\|^2) g(x)] dx
\end{equation}
in the case of $d =2$,  where $\Delta$ is the usual Laplace operator. A 
family of explicit orthonormal basis is constructed in \cite{X06}, which we 
shall follow to construct the basis for the inner products \eqref{ip1} and 
\eqref{ip2}. Both the bases in \cite{X06} and those we shall construct in 
this paper depend on Jacobi  polynomials. It is interesting to note that 
the basis for $\la f, g\ra_{\cI}$ and the one for $\la f, g\ra_\Delta$ have
the same structure and their difference appears in the parameters of the
Jacobi polynomials. 

The paper is organized as follow. In Section 2 we construct explicit 
orthonormal bases for the inner products $\la \cdot, \cdot \ra_{\cI}$  and
 $\la \cdot, \cdot \ra_{\cII}$.  An interesting consequence of the explicit
 formula shows that the orthogonal expansion of a function $f$ with respect
 to these inner products can be computed without involving the derivatives 
 of $f$,  which and the Parseval type relation for the gradient of a function
 will be discussed in Section 3.

\section{Sobolev Orthogonal Polynomials}\label{Sobolev}
\setcounter{equation}{0}

The unit ball in $\RR^d$ is $B^d:=\{x: \|x\| \le 1\}$. Its surface is 
$S^{d-1}:=\{x:  \|x\| =1\}$. Let $d\omega$ denote the Lebesgue measure on 
$S^{d-1}$ and denote the area of $S^{d-1}$ by $\omega_{d}$,  
$$
    \o_{d} := \int_{S^{d-1}} d\o = 2 \pi^{d/2}/\Gamma(d/2).
$$

\subsection{Background and Preliminary}

Let $\Pi^d = \RR[x_1,\ldots,x_d]$ be the ring of polynomials in $d$ variables
and let $\Pi_n^d$ denote the subspace of polynomials of total degree at 
most $n$.  We consider the inner product defined on the polynomial space 
by
\begin{equation}\label{ip1A}
  \langle f,g \rangle_\cI :=  \frac{\l}{\o_d}  \int_{B^d} 
     \na  f(x) \na  g(x) dx + \frac{1}{\o_{d}} \int_{S^{d-1}} f(x)g(x) d\o(x),
\end{equation}
where $\l > 0$, which is the same as \eqref{ip1} up to a normalization; here
it is normalized so that  $\la 1, 1 \ra_{\cI} =1$.  Similarly, we define
\begin{equation}\label{ip2A}
  \langle f,g \rangle_{\cII} := \frac{ \l}{\o_d}  \int_{B^d} 
     \na  f(x) \na  g(x) dx +  f(0)g(0),
\end{equation}
which is the same as \eqref{ip2} up to a normalization chosen here so that 
$\la 1, 1 \ra_{\cII} =1$.  It is easy to see that both inner products are well
defined and positive definite on $\Pi^d$. Let $\la \cdot, \cdot \ra $ denote 
either one of these inner product. 
A polynomial $P$ is orthogonal with respect to $\la \cdot,\cdot \ra$  if it 
is orthogonal to all polynomials of lower degrees with respect to 
$\la \cdot,\cdot \ra$. Denote by $\CV^d(\cI)$ and $\CV^d(\cII)$ the spaces 
of orthogonal polynomials with respect to $\la f,g \ra_\cI$ and 
$\la f,g \ra_\cII$, respectively. Let  $\CV_n^d$ denote either one of these 
spaces. The general theory of orthogonal polynomials in several variables 
(\cite{DX}) shows that the dimension of $\CV_n^d$ is $\binom{n+d-1}{d-1}$.  
Note that the definition puts no restriction on polynomials of the same degree.
A basis $\{P_\alpha \}$ of $\CV_n^d$ is called a mutually orthogonal basis if 
$\la P_\a, P_\beta \ra = 0$ whenever $\a \ne \beta$; it is called an orthonormal 
basis if, in addition, $P_\alpha$ is normalized so that $\la P_\a,P_\a \ra =1$
for all $\a$. 

We can also consider another inner product defined on the unit sphere 
$S^{d-1}$ by 
\begin{equation} \label{ipSd}
  \la f, g \ra_{S}: = \frac{\l}{\o_d} \int_{S^{d-1}} \frac{d }{dr} f(y)
     \frac{d }{dr} g(y) d\o(y) +\frac{1}{\o_d} \int_{S^{d-1}} f(y)g(y) d\o(y)
\end{equation}
where $\lambda \ge 0$ and $d/dr$ denote the derivative in the radial direction.
Here and in the following we use the spherical polar coordinates $x = r x'$ 
for $x \in \RR^d$, $r\ge 0$, and $x' \in S^{d-1}$. 
If $\lambda =0$, then $\la \cdot,\cdot\ra_{S}$ becomes the usual inner product
of $L^2(S^{d-1})$, whose orthogonal polynomials are the spherical harmonics
which we now describe. 

Let $\CH_n^d$ denote the space of homogeneous harmonic polynomials of 
degree $n$, which are homogeneous polynomials of degree $n$ satisfying 
the equation $\Delta P =0$. It is well known that 
$$
\dim \CH_n^d = \binom{n+d-1}{d-1}-\binom{n+d-3}{d-1}:= \s_n.
$$ 
The restriction of $Y \in \CH_n^d$ on $S^{d-1}$ are the spherical harmonics. 
For $Y \in \CH_n^d$ we use the notation $Y(x)$ to denote the harmonic 
polynomials and use $Y(x')$ to denote the spherical harmonics. Since $Y$ is 
a homogeneous polynomial, $Y(x) = r^nY(x')$. Spherical harmonics 
are orthogonal on $S^{d-1}$. Throughout this paper, we use the notation 
$\{Y_\nu^n: 1 \le \nu \le \sigma_n\}$ to denote an orthonormal basis for 
$\CH_n^d$, that is,  
\begin{equation} \label{eq:harmonics}
\frac{1}{\omega_{d}} \int_{S^{d-1}} Y_\mu^n(x') Y_\nu^m(x') d\omega(x')
    = \delta_{\mu,\nu} \delta_{n,m} , \qquad 1 \le \mu,\nu\le \sigma_n.
\end{equation}
In terms of the spherical polar coordinates, $x = r x'$, $r > 0$ and $x' \in S^{d-1}$, 
the Laplace operator can be written as
\begin{equation}\label{eq:Delta}
  \Delta = \frac{\partial^2}{\partial r^2}+ \frac{d-1}{r} 
      \frac{\partial}{\partial r} + \frac{1}{r^2} \Delta_0,
\end{equation}
where $\Delta_0$ is the spherical Laplacian on $S^{d-1}$. 

The Sobolev orthogonal polynomials with respect to the inner product
$\la \cdot,\cdot \ra_S$ defined at \eqref{ipSd} are in fact still spherical 
harmonics. In other words, the space of homogeneous polynomials 
orthogonal with respect to $\la \cdot,\cdot\ra_S$ is exactly $\CH_n^d$. 
More precisely, $\{\sqrt{\l n^2 +1} \,Y_\nu^n: 1 \le \nu \le \s_n\}$ forms an orthonormal 
basis for the space of homogeneous orthogonal polynomials with 
respect to $\la \cdot,\cdot\ra_S$. This follows easily from the fact that
$\frac{d}{dr} Y_\nu^n(x) \vert_{r=1} = n Y_\nu^n(x')$ as $Y_\nu^n$ is 
a homogeneous polynomial of degree $n$. 

It turns out that the Sobolev orthogonal polynomials on the ball have structures 
similar to those of ordinary orthogonal polynomials. Let us consider the 
inner product 
$$
\langle f,g \rangle_\mu : = c_\mu \int_{B^d} f(x) g(x) W_\mu(x)dx,\qquad
   W_\mu(x) = (1-\|x\|^2)^{\mu}, 
$$
where $\mu > -1$ and $c_\mu$ is the normalization constant of $W_\mu$. Let 
$\CV_n^d(W_\mu)$ denote the space of orthogonal polynomials of degree $n$. 
Then a mutually orthogonal basis for $\CV_n(W_\mu)$ is given by (\cite{DX})
\begin{equation}\label{eq:Wmu-basis}
 P_{j,\nu}^n(W_\mu; x) = P_j^{(\mu, n-2j+\frac{d-2}{2})}(2\|x\|^2 -1) 
   Y_{\nu}^{n-2j}(x), \quad 0 \le j \le n/2,
\end{equation}
where $P_j^{(\a,\b)}$ denotes the Jacobi polynomial of degree $j$, which 
is orthogonal with respect to $(1-x)^\alpha(1+x)^\beta$ on $[-1,1]$, and
$\{Y_\nu^{n-2j}: 1 \le j \le \sigma_{n-2j}\}$ is an orthonormal basis for 
$\CH_{n-2j}^d$. 

In \cite{X06} we found a family of Sobolev polynomials with respect to 
$\la \cdot, \cdot \ra_\Delta$ in the form of \eqref{eq:Wmu-basis}. More
precisely, let $\CV_n^d(\Delta)$ denote the space of orthonormal 
polynomials of degree $n$ with respect to the inner product in 
\eqref{innerDelta}, then we proved the following result. 

\begin{thm}
A mutually orthogonal basis $\{Q_{j,\nu}^n: 0 \le j \le \frac{n}{2}, 
1 \le \nu \le \s_{n-2j} \}$for $\CV_n^d(\Delta)$ is given by
\begin{align}\label{eq:basis}
\begin{split}
Q_{0,\nu}^n(x) & = Y_\nu^n(x), \quad \\
 Q_{j,\nu}^n(x) & = (1-\|x\|^2) P_{j-1}^{(2,n-2j+\frac{d-2}{2})}(2\|x\|^2-1)
     Y_\nu^{n-2j}(x), \quad 1 \le j \le \frac{n}{2},
\end{split}
\end{align}
where $\{Y_\nu^{n-2j}: 1 \le \nu \le \sigma_{n-2j}\}$ is an orthonormal basis 
of $\CH_{n-2j}^d$. Furthermore, 
\begin{align}\label{eq:Qnorm}
\langle Q_{0,\nu}^n, Q_{0,\nu}^n \rangle_\Delta = \frac{2n+d}{d}, \qquad
\langle Q_{j,\nu}^n, Q_{j,\nu}^n \rangle_\Delta =\frac{8 j^2(j+1)^2}{d(n+d/2)}.
\end{align}
\end{thm}

\begin{rem} \label{remark1}
A corollary of the theorem shows 
$\CV_n^d(\Delta) = \CH_n^d \oplus (1-\|x\|^2)\CV_{n-2}^d(W_2)$. In
fact, $Q_{j,\nu}(x)$ in \eqref{eq:basis} for $j\ge 1$ can 
be written as
\begin{align}\label{eq:P-basis}
    Q_{j,\nu}^n(x)  = (1-\|x\|^2) P_{j-1,\nu}^{n-2}(W_2;x), \quad 1 \le j \le \frac{n}{2}. 
\end{align}
\end{rem}

The proof of this result relies on the action of $\Delta$ on $Q_{j,\nu}$. More
generally, let 
\begin{equation}\label{eq:Q-basis}
R_{j,\nu}^n(x) := q_j(2\|x\|^2-1) Y_\nu^{n-2j}(x), \qquad 0 \le 2j \le n, 
     \quad Y_\nu^{n-2j} \in \CH_{n-2j}^d,
\end{equation} 
where $q_j$ is a polynomial of degree $j$ in one variable; then the
following lemma holds (see \cite[Lemma 2.1]{X06}). 

\begin{lem} \label{lem:2.1}
Let $R_{j,\nu}^n$ be defined as above. Then
$$
 \Delta\left[(1-\|x\|^2)R_{j,\nu}^n(x) \right] =
    4 \left(\CJ_\beta q_j\right)(2r^2-1) Y_\nu^{n-2j}(x), 
$$   
where $\beta = n-2j + \frac{d-2}{2}$ and
$$
(\CJ_\beta q_j)(s) = (1-s^2) q_j''(s) + (\beta -1 - (\beta+3) s )q_j'(s) 
     - (\beta+1) q_j(s). 
$$
\end{lem}

\subsection{Sobolev orthogonal polynomials with respect to 
$\la \cdot,\cdot\ra_{\cI}$} 

The main result in this subsection is a family of mutually orthogonal basis 
for the inner product in \eqref{ip1A}. 

 \begin{thm}
A mutually orthogonal basis $\{U_{j,\nu}^n: 0 \le j \le \frac{n}{2}, 
1 \le \nu \le \s_{n-2j} \}$ for $\CV_n^d(\cI)$ is given by
\begin{align}\label{eq:basisI}
\begin{split}
 U_{0,\nu}^n(x) & = Y_\nu^n(x), \quad \\
 U_{j,\nu}^n(x) & = (1-\|x\|^2) P_{j-1}^{(1,n-2j+\frac{d-2}{2})}(2\|x\|^2-1)
     Y_\nu^{n-2j}(x), \quad 1 \le j \le \frac{n}{2},
\end{split}
\end{align}
where $\{Y_\nu^{n-2j}: 1 \le \nu \le \sigma_{n-2j}\}$ is an orthonormal basis 
of $\CH_{n-2j}^d$. Furthermore, 
\begin{align}\label{eq:RnormI}
\langle U_{0,\nu}^n, U_{0,\nu}^n \rangle_\cI = n \l +1, \qquad
\langle U_{j, \nu}^n, U_{j, \nu}^n \rangle_\cI =\frac{2 j^2}{n+\frac{d-2}{2}} \l.
\end{align}
\end{thm}

\begin{proof}
A standard argument as in the case for $P_{j,\nu}^n$ shows that 
$\{U_{j,\nu}^n: 0 \le j \le n/2,1 \le \mu \le \s_{n-1}\}$ is a basis for $\Pi_n^d$. 
In order to establish the orthogonality we start with Green's identity, 
$$
   \int_{B^d} \nabla f(x) \nabla g(x) dx = \int_{S^{d-1}}f(x) \frac{d}{dr} g(x) d\o
        - \int_{B^d} f(x) \Delta g(x) dx,
$$
where $d/dr$ is the normal derivative which coincides with the derivative
in the radial direction. This identity can be used to rewrite the inner product 
$\la \cdot, \cdot \ra_\cI$ as  
\begin{equation} \label{eq:ip1B}
\la f, g \ra_{\cI} = \frac{1}{\o_d} \int_{S^{d-1}} f(x) 
         \left[ \l \frac{d}{dr}g(x) + g(x)\right] d\o 
               -   \frac{\l}{\o_d} \int_{B^d} f(x) \Delta g(x) dx.
\end{equation}
First we consider the case $j=0$; that is, the orthogonality of $U_{0,\nu}^n
 = Y_\nu^n$. Setting $j=0$ in \eqref{eq:Wmu-basis} shows that $Y_\nu^n$ 
is an orthogonal polynomial in $\CV_n^d(W_\mu)$.  Since 
$U_{j,\nu}^n(x) \vert_{r=1} =0$, $\frac{d}{dr} U_{j,\nu}^n(x) \vert_{r=1} 
    = -2 P_{j-1}^{(1,n-2j+\frac{d-2}{2})}(1)$ and $\Delta U_{j,\nu}^m \in 
 \Pi_{m-2}^d$,  it follows from \eqref{eq:ip1B}
that for $m < n$, $j \ge 0$ and $0 \le \mu \le \s_{m-2j}$, 
$$
 \la U_{0,\nu}^n, U_{j,\mu}^m\ra_{\cI} = -2 P_{j-1}^{(1,n-2j+\frac{d-2}{2})}(1)
  \frac{1}{\o_{d}} \int_{S^{d-1}} Y_\nu^n(x') Y_\mu^m(x') d\o(x') =0.
$$
Furthermore, using the fact that $\frac{d}{dr} Y_\nu^n(x)\vert_{r=1} = 
nY_\nu^n(x')$, the same consideration shows that 
$$
  \la U_{0,\nu}^n, U_{0,\nu}^n\ra_{\cI} =  (\l n+1)  
      \frac{1}{\o_{d}} \int_{S^{d-1}} \left[Y_\nu^n(x')\right]^2 d\o(x') =\l n+1.
$$
Next we consider $U_{j,\nu}^n$ for $j \ge 1$. In this case $U_{j,\nu}(x)\vert_{r=1} 
=0$ since it contains the factor $(1-\|x\|^2)$, which is zero on $S^{d-1}$. 
Consequently, the first term in
\begin{align*}
 \la U_{j,\nu}^n, U_{l,\mu}^m\ra_{\cI} = & \frac{1}{\o_d}\int_{S^{d-1}}
   U_{j,\nu}^n(x') \left[ \frac{d}{dr} U_{l,\nu}^m+ U_{l,\nu}^m\right](x') d\o(x') \\
   &- \frac{\l}{\o_d} \int_{B^d} U_{j,\nu}^n(x) \Delta U_{l,\mu}^m(x) dx
\end{align*}
is zero. For the second term, we use Lemma \ref{lem:2.1} to derive a 
formula for $\Delta U_{j,\nu}^n$. The formula in the lemma gives 
$$
 \Delta U_{j,\nu}^n(x) = 4 \left (\CJ_\b P_{j-1}^{(1,\b)}\right)(2r ^2 -1)
       Y_\nu^{n-2j}(x), \qquad \b = n-2j+\frac{d-2}{2}. 
$$
On the other hand, the Jacobi polynomial $P_{j-1}^{(1,\b)}$  satisfies 
the differential equation (\cite[p. 60]{Sz})
$$
     (1-s^2) y''-(1-\b+(3+\b)s)y'+(j-1)(j+\b+1)y=0, 
$$
which implies that $(\CJ_\b P_{j-1}^{(1,\b)})(s)= -j (j+\b)P_{j-1}^{(1,\b)}(s)$.
Consequently, we obtain
\begin{align} \label{eq:DR=P}
   \Delta U_{j,\nu}^n(x) = & \, - 4j(j+\b)P_{j-1}^{(1,\b)}(2r^2-1)Y_\nu^{n-2j}(x) \\
         = &\, - 4j(n-j+\tfrac{d-2}{2})P_{j-1,\nu}^{n-2}(W_1;x). \notag
\end{align}
Hence, using the fact that $U_{l,\mu}^m(x) = (1-\|x\|^2) P_{l-1,\mu}^{m-2}(W_1;x)$
(see \eqref{eq:R-P} below), we derive from \eqref{eq:DR=P} that 
\begin{align*}
  & \int_{B^d} U_{l,\mu}^m(x) \Delta U_{j,\nu}^n(x) dx \\
  & \qquad = 
     -4 j (n-j+\tfrac{d-2}{2}) \int_{B^d} P_{l-1,\mu}^{m-2}(W_1;x) 
        P_{j-1,\nu}^{n-2}(W_1;x) (1-\|x\|^2)dx \\
  & \qquad = -4 j (n-j+\tfrac{d-2}{2}) \int_{B^d} \left[P_{j-1,\nu}^{n-2}(W_1;x) \right]^2
       (1-\|x\|^2)dx\, \delta_{n,m} \delta_{j,l} \delta_{\nu,\mu}. 
\end{align*}
The norm of $P_{j,\nu}^n(W_\mu;x)$ in $L^2(B^d,W_\mu)$ can be computed 
as in \cite[p. 39]{DX}, where the weight function is normalized. The result is 
\begin{equation}\label{eq:W1-norm}
       \frac{1}{\o_d} \int_{B^d} \left[P_{j-1,\nu}^{n-2}(W_1;x) \right]^2 (1-\|x\|^2)dx 
        = \frac{1}{2} \frac{j}{(n-j+\frac{d-2}{2})(n+\frac{d-2}{2})}. 
\end{equation}
Putting these pieces together proves \eqref{eq:RnormI} for $j \ge 1$ and the 
theorem.
\end{proof} 

From the explicit formula of the basis \eqref{eq:basisI} it follows that 
$U_{j,\nu}^n$ is related to orthogonal polynomials with respect to
$W_1(x) = 1-\|x\|$. In fact, we have
\begin{equation} \label{eq:R-P}
       U_{j,\nu}^n(x) = (1-\|x\|^2) P_{j-1,\nu}^{n-2} (x), \qquad j \ge 1,
\end{equation}
which has already been used in the above proof. An immediate consequence
is the following corollary.

\begin{cor}
For $n \ge 1$, 
$$
  \CV_n^d(\cI) = \CH_n^d \oplus (1-\|x\|)\CV_{n-2}^d(W_1).
$$
\end{cor}

\begin{rem} 
It is interesting to compare this relation to the analogous relation between 
$Q_{j,\nu}^n \in \CV_n^d(\Delta)$ and orthogonal polynomials with respect 
to $W_2(x)$ given in \eqref{eq:P-basis}. See Remark \ref{remark1}. 
\end{rem}

It is known that polynomials in $\CV_n^d(W_1)$ are eigenfunctions of
a second order differential operator (cf. \cite{DX}), more precisely, 
$$
 \CD P: = \left(\Delta- \la  x, \nabla \ra^2 - (d+1)\la x, \nabla \ra \right) P
            = -(n+d)(n+2) P
$$
for $P \in \CV_n^d (W_1)$. Using the relation \eqref{eq:DR=P}, the 
following result is immediate. 

\begin{cor}
Polynomials in $\CV_n^d(\cI)$ satisfy the differential equation
$$
      [\CD  + (n+d)(n+2)]\Delta P = 0. 
$$
\end{cor}

\subsection{Sobolev orthogonal polynomials 
with respect to $\la \cdot,\cdot\ra_{\cII}$} 

We turn our attention to the inner product in \eqref{ip2A}. Again our 
main result is an explicit family of mutually orthogonal basis.  The basis that 
we will give for $\CV_n^d(\cII)$ turns out to be similar to the basis for 
$\CV_n^d(\cI)$ given in \eqref{eq:basisI}. In fact, for $n$ is odd, the two bases 
are identical, whereas for $n$ is even, the two bases differ by just one 
element. We will need a result on Sobolev orthogonal polynomials of 
one variable with respect to the inner product 
$$
 (f,g) :=  2^{2-d/2} \l \int_{-1}^1  f'(s) g'(s) (1+s)^{d/2} ds + f(-1)g(-1).  
$$
Let us denote by $q_k$ the $k$-th orthogonal polynomial of one variable with 
respect to the above inner product. It is easy to see (cf. \cite{JKL}) that 
$$
  q_0(x) =1, \qquad q_k(x) = \int_{-1}^x P_{k-1}^{(0,\frac{d}{2})}(t)dt, \qquad k \ge 1,
$$
is an orthogonal basis with respect to $(\cdot,\cdot)$. The Jacobi polynomials 
$P_k^{(-1,\frac{d-2}{2})}(x)$ are well defined for $k \ge 1$ and we 
have the relation \cite[(4.5.5), p. 72]{Sz}
$$
\frac{d}{dx} P_k^{(-1,\frac{d-2}{2})}(x) = \frac{1}{2}\left(k+ \frac{d-2}{2} \right) 
       P_{k-1}^{(0,\frac{d}{2})}(x),
$$
from which we immediately deduce that 
\begin{equation} \label{qk}
  q_k(x) = \frac{2}{k+\frac{d-2}{2}} \left(P_k^{(-1,\frac{d-2}{2})}(x)
             - (-1)^k \frac{(d/2)_k}{k!}\right), \qquad k \ge 1,
\end{equation}
where $(a)_k=a (a+1) \ldots (a+k-1)$ is the shifted factorial and we have used
the fact that $P_k^{(-1,\b)}(-1) = (-1)^k (\b+1)_k/k!$. 

\begin{thm}
A mutually orthogonal basis $\{V_{j,\nu}^n: 0 \le j \le \frac{n}{2}, 1 \le \nu \le \s_{n-2j} \}$
for $\CV_n^d(\cII)$ is given by
\begin{align}\label{eq:basisII}
\begin{split}
 V_{j,\nu}^n(x) & = U_{j,\nu}^n(x), \quad 1 \le j \le \lfloor \tfrac{n-1}{2} \rfloor,\\
 V_{\frac{n}{2}}^n(x) & =  \frac{1}{n+\frac{d-2}{2}} 
        \left(P_{\frac{n}{2}}^{(-1,\frac{d-2}{2})}(2\|x\|^2-1) 
         -(-1)^{\frac{n}{2}} \frac{(d/2)_{\frac{n}{2}}} {(n/2)!}
         \right), 
\end{split}
\end{align}
where $V_{\frac{n}{2}}^n(x) := V_{\frac{n}{2},\nu}^n(x)$ holds only when $n$ 
is even. Furthermore,
\begin{align}\label{eq:RnormII}
\begin{split}
\la V_{j,\nu}^n, V_{j,\nu}^n \rangle_\cII &\  = \la U_{j,\nu}^n, U_{j,\nu}^n \rangle_\cI,
      \quad 1 \le j \le \lfloor \tfrac{n-1}{2} \rfloor, \\ 
\la V_{\frac{n}{2}}^n, V_{\frac{n}{2}}^n \rangle_\cII & \ =\frac{8 \l }{n  + \frac{d-2}{2}}.
\end{split}
\end{align}
\end{thm}

\begin{proof}
First we note that, for $n$ even, there is only one element in $V_{\frac{n}{2},\nu}^n$,
which is why we drop the $\nu$ in the notation. For $1 \le j \le (n-1)/2$,  the explicit expression of the polynomials $U_{j,\nu}$  at \eqref{eq:basisI} shows that 
$$
 \int_{S^{d-1}} V_{j,\nu}^n(x) V_{l,\mu}^m(x) d\o = 0 =   
       V_{j,\nu}^n(0) V_{l,\mu}^m(0)
$$
for $0 \le l \le m/2$ and $m \le n$, where the left hand side is zero because
the presence of $(1-\|x\|^2)$, whereas the right hand side is zero because 
$Y_\nu^{n-2j}(0) =0$ for $j < n/2$. Thus, $\la V_{j,\nu}^n, V_{l,\mu}^m \rangle_\cII
 = \la U_{j,\nu}^n, U_{l,\mu}^m \rangle_\cI$ for $0 \le l \le m/2$ and $m \le n$.
Hence, we only need to consider the orthogonality of $V_{\frac{n}{2}}^n$. 
 
First assume that $2 j \ne m$. By Green's formula and \eqref{eq:DR=P}, 
\begin{align*}
  \la V_{\frac{n}{2}}^n, V_{j,\nu}^m\ra_{\cII} = &\ 
          V_{\frac{n}{2}}^n(0) V_{j,\nu}^m(0)+ \l  V_{\frac{n}{2}}^n(1)
    \frac{1}{\o_d} \int_{S^{d-1}} \frac{d}{dr} V_{j,\nu}^m(x) d\o \\
  & + \l 4 j(m-j+ \tfrac{d-1}{2}) \int_{B^d} V_{\frac{n}{2}}^n(x) 
        P_{j-1}^{(1,m-j+ \frac{d-1}{2})}(2\|x\|^2-1)Y_\nu^{m-2j}(x) dx 
\end{align*}
The first term is zero as before. The second term is zero by the orthogonality 
of $Y_\nu^{m-2j}$ since $\frac{d}{dr} V_{j,\nu}^m(x) \vert_{r=1} =  c Y_\nu^{m-2j}(x')$ 
and $m- 2j \ne 0$. The third term is also zero,  since $V_{\frac{n}{2}}^n$ is a 
radial function so that we can use the formula 
$$
    \int_{B^d} g(x) d x = \int_0^1 \int_{S^{d-1}} g (r x') d\o(x') r^{d-1} dr
$$
and the orthogonality of $Y_\nu^{n-2j}$. Next we consider the case of 
$m = 2j$. By definition, $V_{\frac{n}{2}}^n = q_{\frac{n}{2}}(2\|x\|^2-1)$ with
$q_\k$ given in \eqref{qk}. We use the fact that for $f, g: \RR \mapsto \RR$, 
$$
       \nabla f(2 \|x\|^2 -1)\nabla g(2 \|x\|^2-1) = 16 \|x\|^2
         f'(2 \|x\|^2 -1) g'(2 \|x\|^2-1), 
$$ 
and the definition of $\la \cdot,\cdot \ra_{\cII}$ in \eqref{ip2A} to conclude that 
\begin{align*}
\la V_{\frac{n}{2}}^n, V_{\frac{m}{2}}^m\ra_{\cII} & = 16 \l \int_0^1 r^{d+1}
    q_{\frac{n}{2}}'(2r^2 -1)q_{\frac{m}{2}}'(2r^2 -1) dr 
     + q_{\frac{n}{2}} (-1)q_{\frac{m}{2}}(-1) \\
& = 2^{2-d/2} \l \int_{-1}^1 q_{\frac{n}{2}}'(t)q_{\frac{m}{2}}'(t) (1+t)^{d/2} dt
  + q_{\frac{n}{2}} (-1)q_{\frac{m}{2}}(-1), 
\end{align*}
which is zero whenever $m \ne n$ by the definition of $q_k$. In the case of 
$n =m \ge 1$, we have 
\begin{align*}
\la V_{\frac{n}{2}}^n, V_{\frac{m}{2}}^m\ra_{\cII}  = 
 2^{2-d/2} \l \int_{-1}^1 \left[P_{n-1}^{(0,\frac{d}{2})}(t) \right]^2  (1+t)^{d/2} dt 
  = \frac{8\l}{n + \frac{d-2}{2}},
\end{align*}
using the well-known norm of the Jacobi polynomials (\cite[p. 68]{Sz}). 
\end{proof}


\section{Expansions in Sobolev Orthogonal Polynomials}\label{Expand}
\setcounter{equation}{0}

Let $\CV_n^d$ be the space of orthogonal polynomials of degree $n$ with 
respect to an inner product $\la \cdot,\cdot \ra$. Let $H^2(B^d)$ denote the 
space of functions for which $\langle f,f \rangle$ is finite.  In the case of 
$\CV_n^d(W_\mu)$ with respect to $\la \cdot, \cdot \ra_\mu$, the space
is just $L^2(B^d, W_\mu)$. In the case of either $\la \cdot,\cdot\ra_\cI$ 
or $\la \cdot,\cdot \ra_\cII$, the presence of the derivative in the definition of 
the inner product shows that $H^2(B^d)$ is not the $L^2$ space on $B^d$.
Nevertheless , the standard Hilbert space theory shows that every $f\in H^2(B^d)$ 
can be expanded into a series in Sobolev orthogonal polynomials. Thus, in
all these cases we have
$$
        H^2(B^d) = \sum_{n=0}^\infty\oplus \CV_n^d: 
                           \qquad f = \sum_{n=0}^\infty \proj_{\CV_n^d} f,
$$
where $\proj_{\CV_n^d}: H^2(B^d) \mapsto \CV_n^d$ is the projection operator.
Let $\{P_{j,\nu}^n\}$ denote a mutually orthogonal basis for $\CV_n^d$, such
as one of the bases \eqref{eq:Wmu-basis}, \eqref{eq:basisI} or \eqref{eq:basisII}. 
The projection operator satisfies 
\begin{equation} \label{eq:proj}
 \proj_{\CV_n^d} f(x)=\sum_{0 \le j \le n/2}  \sum_{\nu =0}^{\sigma_{n-2j}}
  H_{j,\nu}^{-1}    \wh f_{j,\nu}^n  P_{j,\nu}^n(x), \qquad \wh f_{j,\nu}^n = 
            \la f,P_{j,\nu}^n \ra_\cI,
\end{equation}
where $H_{j,\nu} = \langle P_{j,\nu}^n, P_{j,\nu}^n \rangle_\cI$. Let $P_n(x,y)$ 
denote the reproducing kernel of $\CV_n^d$. In terms of the orthonormal
basis, the reproducing kernel can be
written as 
$$
  P_n(x,y) = \sum_{0\le j \le n/2} H_j^{-1} 
             \sum_\nu U_{j,\nu}^n(x) U_{j,\nu}^n(y).  
$$
The projection operator can be written as an integral operator with $P_n(\cdot,\cdot)$ 
as its kernel, which means that 
$$
    \proj_{\CV_n^d} f = \la f,  P_n(x, \cdot)\ra. 
$$
To distinguish between various inner products, we shall use the notation such 
as $\proj_{\CV_n^d(W_\mu)} f$ and $\proj_{\CV_n^d(\cI)} f$ to denote the 
corresponding projection operator.

We consider the inner product $\la \cdot,\cdot \ra_\cI$ first. In this case, 
we denote the reproducing kernel by $P_n^{\cI}(\cdot,\cdot)$ and we have 
\begin{align*}
   \proj_{\CV_n^d(\cI)}  f(x)
    = \frac{\l}{\o_d}  \int_{B^d}\nabla f(y) \nabla P_n^\cI(x,y) dy 
       + \frac{1}{\o_d} \int_{S^{d-1}} f(y) P^\cI_n(x,y) d\o(y),
\end{align*}
where $\nabla$ is applied on $y$ variable. Note also that $H_{j\nu} = H_j$ since 
they are independent of $\nu$ as shown in \eqref{eq:RnormI}. 

Our first result shows that the orthogonal expansion can be computed 
without involving derivatives of $f$. Recall that $\{Y_\nu^n\}_\nu$ denotes 
an orthonormal basis for $\CH_n^d$. We denote the Fourier coefficients 
of $f \in L^2(S^{d-1})$ with respect to this basis of spherical harmonics by
$$
\la f,Y_\nu^n\ra_{L^2(S^{d-1})}:=\frac{1}{\omega_d}\int_{S^{d-1}} 
    Y_\nu^{n}(y')f(y') d\omega(y').
$$
Furthermore, recall that $P_{j,\nu}^n(W_\mu)$ given in \eqref{eq:Wmu-basis} 
forms a mutually orthogonal basis for $\CV_n^d(W_\mu)$. 

\begin{prop} \label{FourierCoeff}
The Fourier coefficients for $\la \cdot,\cdot \ra_\cI$ satisfy, 
\begin{align} \label{eq:FourierI0}
 \wh f_{0,\nu}^n = (\l n +1)\la f,Y_\nu^n\ra_{L^2(S^{d-1})}, \quad 
 \end{align} 
and for $j > 0$, 
\begin{align} \label{eq:FourierI}
 \wh f_{j,\nu}^n =  &\, - 2 j \l \la f, Y_\nu^{n-2j} \ra_{L^2(S^{d-1})}\\
 &   + 4j \left(n-j+ \frac{d-2}{2}\right) \frac{\l}{\o_d} 
         \int_{B^d} f(y) P_{j-1,\nu}^{n-2} (W_1;y) dy.
\notag
\end{align} 
\end{prop}

\begin{proof}
From Green's identity we immediately obtain that 
$$
  \la f ,U_{j,\nu}^n \ra_\cI =  -\frac{\l}{\o_d} \int_{B^d} f(x) \Delta U_{j,\nu}^n(x)dx
       + \frac{1}{\o_d} \int_{S^{d-1}}
             f(x) \left[ \l \frac{d}{dr} U_{j,\nu}^n(x) + U_{j,\nu}^n(x)\right] d\o. 
$$
Using the explicit formula of $U_{j,\nu}^n$ at \eqref{eq:basisI} and the fact
that $Y_\nu^n$ is a homogeneous polynomial of degree $n$, it is easy to
see that 
$$
   \frac{d}{dr} U_{0,\nu}^n(x) \Big\vert_{r=1} =  n Y_\nu^n(x'), \quad 
   \frac{d}{dr} U_{j,\nu}^n(x) \Big\vert_{r=1} =  -2j  Y_\nu^{n-2j}(x'), \quad j \ge 1, 
$$
where we have used the fact that $P_{j-1}^{(1,\beta)}(1) = j$. Furthermore,
$U_{0,\nu}^n(x)\vert_{r=1} = Y_\nu^n(x')$ and $U_{j,\nu}^n(x)\vert_{r=1} =0$
for $j \ge 1$. Thus, 
$$
   \wh f_{0,\nu}^n  = \la f, U_{0,\nu}^n\ra_\cI = 
       (\l n+1) \frac{1}{\o_d} \int_{S^{d-1}} f(y)Y_\nu^n(y) d\o
$$
using the fact that $\Delta Y_\nu^n =0$. For $j \ge 1$,  the stated result follows 
from the relations derived above,  equations\eqref{eq:R-P} and \eqref{eq:DR=P}.
\end{proof} 

Let us denote by $P_n(W_\mu; x,y)$ the reproducing kernel of $\CV_n^d(W_\mu)$,
which can be written as 
$$
 P_n(W_\mu; x,y) = \sum_{|\alpha| = n} A_{\alpha,\mu}^{-1} 
   P_\alpha(W_\mu;x) P_\alpha(W_\mu;y), 
$$
where $A_{\alpha,\mu} = c_\mu \int_{B^d}[P_\alpha(W_\mu;y)]^2 W_\mu(y)dy$ in
which $c_\mu$ is the normalization of $W_\mu$. A compact formula for this
kernel can be found in \cite{X99}.  We also denote the projection operator
of $f \in L^2(S^{d-1})$ onto $\CH_n^d$ by $Y_n f$.  It is well known that 
$$
Y_n f(x) = \|x\|^n \frac{n+\frac{d-2}{2}} {\frac{d-2}{2}} \frac{1}{\o_d}
   \int_{S^{d-1}} f(y) C_n^{\frac{d-2}{2}} ( x' \cdot y ) d\o(y), 
$$
for $x \in B^d$ and $x'= x/\|x\| \in S^{d-1}$, where $C_n^\lambda(t)$ denote 
the Gegenbauer polynomial of degree $n$ and $x \cdot y$ is the usual dot 
product of $x, y \in \RR^d$. 

\begin{thm} \label{thm:proj}
For $f \in H^2(B^d)$ and $x \in B^d$,
\begin{align*}
 \proj_{\CV_n^d(\cI)} f(x)& = Y_n f(x) + 
  (1-\|x\|^2) \left[ \frac{d(\frac{d}{2}+1)}{\o_d}  
   \int_{B^d} f(y) P_{n-2}(W_1;x,y) dy  \right. \\
    &  \left. \quad  - \left(n+\frac{d-2}{2} \right)
    \sum_{1 \le j \le n/2} \frac{1} {j}  P_{j-1}^{(1,n-2j+\frac{d-2}{2})} (2\|x\|^2-1)
            Y_{n-2j}f(x) \right]. 
\end{align*} 
\end{thm} 

\begin{proof}
The values of $H_j = \langle U_{j,\nu}^n, U_{j,\nu}^n \rangle_\cI$ are 
given in \eqref{eq:RnormI}, from which and \eqref{eq:FourierI0} it follows
that 
$$
  \sum_{\nu=1}^{\sigma_n} H_0^{-1} \wh f_{0,\nu}^n U_{0,\nu}^n(x)  
  = \frac{1}{\omega_d} \int_{S^{d-1}} f(y') \sum_{\nu=1}^{\sigma_n} 
     Y_\nu^n(y')Y_\nu^n(x) d\omega(y') = Y_n f(x), 
$$
where the last step follows from the summation formula of spherical 
harmonics,
\begin{equation}\label{eq:zono}
 \sum_{\nu=1}^{\sigma_n} Y_\nu^n(x)Y_\nu^n(y)=
   \|x\|^n \sum_{\nu=1}^{\sigma_n} Y_\nu^n(x')Y_\nu^n(y)
   = \|x\|^n \frac{n+\frac{d-2}{2}}{\frac{d-2}{2}}
     C_n^{\frac{d-2}{2}} (x' \cdot y)
\end{equation}
for $x', y \in S^{d-1}$. For $j \ge 1$, it follows from \eqref{eq:FourierI} and 
\eqref{eq:RnormI} that 
\begin{align*}
H_j^{-1} \wh f_{j,\nu}^n 
= &\, -\left( n+\frac{d-2}{2} \right) j^{-1} 
      \frac{1}{\o_d} \int_{S^{d-1}} f(y')Y_\nu^{n-2j}(y') d\omega(y')\\
 &\, +2\left( n+\frac{d-2}{2} \right) \left(n-j+ \frac{d-2}{2} \right) j^{-1} 
  \frac{1}{\o_d} \int_{B^d} f(y) P_{j-1,\nu}^{n-2}(W_1;y)dy,
\end{align*}
Hence, by the formula of $U_{j,\nu}^n$ at \eqref{eq:basisI} and
\eqref{eq:R-P} as well as the identity \eqref{eq:zono}, we obtain that 
\begin{align*}
 \proj_{\CV_n^d(\cI)} f (x) &\, =  Y_n f(x) \\
      &\, -\left( n+\frac{d-2}{2} \right) (1-\|x\|^2)\sum_{1\le j \le n/2} 
          \frac{1}{j} P_{j-1}^{(1,n-2j+\frac{d-2}{2})}(2\|x\|^2-1)Y_{n-2j} f \\
 &\, +2\left( n+\frac{d-2}{2} \right) (1-\|x\|^2) \sum_{1 \le j \le n/2}
      \sum _\nu \left(n-j+ \frac{d-2}{2} \right) j^{-1} \\
  &\qquad\quad \times    \frac{1}{\o_d} \int_{B^d} f(y) P_{j-1,\nu}^{n-2}(W_1;y) dy
          P_{j-1,\nu}^{n-2}(W_1;x).  
\end{align*}
From \eqref{eq:W1-norm} and the fact that the normalization constant of the 
weight function $W_1$ on $B^d$ is given by $d (d/2+1)/\o_d$, we obtain the
norm of $P_{j-1,\nu}^{n-2}(W_1)$ in $L^2 (B^d;W_1)$, from which the stated result
follows from the definition of the reproducing kernel of $\CV_n^d(W_1)$. 
\end{proof}

Despite the presence of $P_{n-2}(W_1;x,y)$, we note that the corresponding
integral in $\proj_{\CV_n^d(\cI)}$ is not the projection of $f $ on $\CV_n^d(W_1)$ 
since the weight function $W_1$ is not in the integral. 

It follows from Theorem \ref{thm:proj} that the orthogonal expansion of $f$ 
with respect to $\langle \cdot,\cdot \rangle_\cI$ coincides with the spherical 
harmonic expansion of $f$ when restricted on $S^{d-1}$. More precisely,
we have 
$$
      \proj_{\CV_n^d(\cI)} f(x)  = Y_n f(x), \qquad x \in  S^{d-1}. 
$$
It is worthwhile to mention that such a result also holds for the projection 
operator with respect to $\la \cdot,\cdot \ra_{\Delta}$ (\cite{X06}).  More 
interesting, however, is the following result: 

\begin{cor}
If $ f(x) = (1-\|x\|^2) g(x) \in H^2(B^d)$, then 
\begin{align}\label{I=W}
    \proj_{\CV_n^d(\cI)} f(x) = (1-\|x\|^2) \proj_{\CV_n^d(W_1)} g(x). 
\end{align}
\end{cor}

\begin{proof}
This follows from the fact that $f(x) =1$ on the sphere, so that $Y_{n-2j} f =0$
for $0 \le j \le n/2$, and the fact that 
$$
  \proj_{\CV_n^d(W_1)} g(x) = 
      \frac{d(\frac{d}{2}+1)}{\o_d}  \int_{B^d} g(y) P_{n-2}(W_1;x,y)(1-\|y\|^2) dy, 
$$
since $d(d/2+1)/\o_d$ is the normalization constant of $W_1(y) = (1-\|y\|^2)$
on $B^d$.
\end{proof}

Recall that we start with the problem of finding a Parseval relation for the 
gradient of a function. Since the Parseval identity holds in $H^2$, it follows 
readily from Proposition \ref{FourierCoeff} and \eqref{eq:RnormI} that 
\begin{align*} 
\la f, f\ra_\cI &\ = \sum_{n=0}^\infty \sum_{ j=0}^{\lfloor \frac{n}{2}\rfloor}
                 \sum_{\nu=1}^{\sigma_{n-2j}} \frac{ \big|\la f,U_{j,\nu}^n\ra_{\cI}\big|^2 }   
                    { \la U_{j,\nu}^n,U_{j,\nu}^n\ra_{\cI}} 
           = \sum_{n=0}^\infty \sum_{ j=0}^{\lfloor \frac{n}{2}\rfloor}
                 \sum_{\nu=1}^{\sigma_{n-2j}} \frac{\big|\wh f_{j,\nu}^n\big|^2 }   
                    { \la U_{j,\nu}^n,U_{j,\nu}^n\ra_{\cI}} \\
             &\ = \sum_{n=0}^\infty (\l n +1) \sum_{\nu=1}^{\sigma_n} 
                 \left |\la f, Y_\nu^n\ra_{L^2(S^{d-1})} \right|^2
            + 2 \sum_{n=0}^\infty  \sum_{j=1}^{\lfloor \frac{n}{2} \rfloor}
                  \left(n+\frac{d-2}{2}\right) \l \\
  &    \qquad   \times  \sum_{\nu=1}^{\sigma_{n-2j}} 
         \left|\la f, Y_\nu^{n-2j} \ra_{L^2(S^{d-1})} - \left(n-j+ \frac{d-2}{2}\right) 
             \la f, P_{j-1,\nu}^{n-2}(W_1) \ra_{L^2(B^d)} \right|^2, 
 \end{align*}
where 
$$
  \la f, P_{j-1,\nu}^{n-2}(W_1)\ra_{L^2(B^d)}: =  
        \frac{2}{\o_d} 
             \int_{B^d} f(y) P_{j-1,\nu}^{n-2} (W_1;y) dy. 
$$
Recall the definition of $\la\cdot,\cdot\ra_\cI$ at \eqref{ip1A}. Dividing the 
above equation by $\l$ and taking the limit $\l \to \infty$, we end
up with the following result: 

\begin{cor} \label{Parseval}
For $f \in H^2$ we have 
\begin{align*} 
 &\frac{1}{\o_d} \int_{B^d} \left|\na f(x)\right|^2 dx  
      = \sum_{n=0}^\infty n \sum_{\nu=1}^{\sigma_n} 
                 \left |\la f, Y_\nu^n\ra_{L^2(S^{d-1})} \right|^2  + 2 \sum_{n=0}^\infty
                   \left(n+\frac{d-2}2\right) 
  \\&\quad   \times  \sum_{j=1}^{\lfloor \frac{n}{2} \rfloor} 
       \sum_{\nu=1}^{\sigma_{n-2j}} 
      \left|\la f, Y_\nu^{n-2j} \ra_{L^2(S^{d-1})} -   \left(n-j+ \frac{d-2}{2}\right) 
           \la f, P_{j-1,\nu}^{n-2}(W_1) \ra_{L^2(B^d)} \right|^2.
 \end{align*}
\end{cor}

This is the Parseval type relation for the gradient of $f$ on $B^d$. The
multiplier $n$ and $2n+d$ shows the impact of the derivative in the 
right hand side. In fact, it is  well known that 
$$
  \frac{1}{\o_d} \int_{S^{d-1}} |f(x)|^2 d\o(x) = 
   \sum_{n=0}^\infty \sum_{\nu=1}^{\sigma_n} 
                 \left |\la f, Y_\nu^n\ra_{L^2(S^{d-1})} \right|^2. 
$$
 
Because of the relation \eqref{I=W} between $\proj_{\CV_n^d(\cI)}f$ and 
$\proj_{\CV_n^d(W_1)}g$ for $f (x) = (1-\|x\|^2) g(x)$, we can expect a 
relation between the Fourier coefficients of $f$ and those of  $g$. Let 
us denote the Fourier coefficients of $g$ in $L^2(B^d,W_1)$ by 
$$
   \{\wh g_{j,\nu}^n\}_{W_1} : =  \frac{2}{\o_d} 
        \int_{B^d} g(y) P_{j,\nu}^n (W_1,y) (1-\|y\|^2)dy. 
$$

\begin{cor}
For $f(x) =(1-\|x\|^2 g(x) \in H^2$ we have 
\begin{align*} 
 \frac{1}{\o_d} \int_{B^d} \left|\na f(x)\right|^2 dx  
   =  2 \sum_{n=0}^\infty \left(n+\frac{d-2}2\right) \sum_{j=1}^{\lfloor \frac{n}{2} \rfloor} 
     \sum_{\nu=1}^{\sigma_{n-2j}} 
            \left(n-j+ \frac{d-2}{2}\right)^2 \left[  \{\wh g_{j,\nu}^n\}_{W_1} \right]^2.
 \end{align*}
\end{cor}

We can also consider the orthogonal expansions with respect to 
$\la\cdot, \cdot \ra_\cII$. However, since the orthogonal polynomials 
$V_{j,\nu}^n$ in \ref{eq:basisII} differ from $U_{j,\nu}^n$ in only one 
element when $n $ is even, there is little left to be done. Using the fact
that $V_{\frac{n}{2}}^n(0)=0$, Green's formula shows that 
$$
  \la f, V_{\frac{n}{2}}^n \ra_\cII = - \frac{4\l}{\o_d}\int_{S^{d-1}} f(y)d\o(y)
  +    \frac{\l}{\o_d}\int_{B^d} f(x) \Delta V_{\frac{n}{2}}^n(x)dx,
$$
where we have used the fact that $\frac{d}{dr} V_{\frac{n}{2}}^n (x) \vert_{r=1} =4$.
Since $V_{\frac{n}{2}}^n(x) = q_{\frac{n}{2}}(2 \|x\|^2-1)$ is a radial function,
one can compute $\Delta V_{\frac{n}{2}}(x)$ by \eqref{eq:Delta}. We note, 
in passing, that the relation such as \eqref{I=W} no long holds for 
$\proj_{\CV_n^d(\cII)} f$. Furthermore, since the Parseval type relation
in Corollary \ref{Parseval} is obtained by dividing $\lambda$ and then taking 
limit $\l \to \infty$, and $\lim_{\l \to \infty} \la f,g\ra_\cI /\l =  
\lim_{\l \to \infty} \la f,g\ra_\cII /\l$, the orthogonal expansion for 
$\la \cdot,\cdot\ra_\cII$ will lead to exactly the same Parseval relation as
the one in Corollary \ref{Parseval}. 

\bigskip\noindent
{\it Acknowledgment.} The author thanks Dr. Greg Forbes  of QED 
Technologies Inc., Sydney, Australia, for drawing his attention to 
this problem.


\begin{thebibliography}{99} 

\bibitem{A}
       K. Atkinson and O. Hansen,
       Solving the nonlinear Poisson equation on the unit disk,
       \textit{J. Integral Equations Appl.}\textbf{17} (2005), 223--241.


\bibitem{DX}
        C. F. Dunkl and Yuan Xu,
        \textit{Orthogonal polynomials of several variables},
        Cambridge Univ. Press, 2001. 
 
\bibitem{GF}
       G. Forbes,   
        Personal communication. 
        
\bibitem{G}
        W. Gautschi, 
        \textit{Orthogonal Polynomials: Computation and Approximation},
        Oxford Univ. Press, 2004.

\bibitem{JKL}
        I. H. Jung, K. H.  Kwon and J. K.  Lee, 
        Sobolev orthogonal polynomials relative to 
         $\lambda p(c)q(c)+\langle\tau,p'(x)q'(x)\rangle$. 
        \textit{Commun. Korean Math. Soc.} \textbf{12} (1997), 603--617.

\bibitem{Sz}
        G. Szeg\H{o},
        \textit{Orthogonal Polynomials},
        Amer. Math. Soc. Colloq. Publ. Vol.23, Providence, 4th edition,
        1975.

\bibitem{X99}
        Yuan Xu,
        Summability of Fourier orthogonal series for Jacobi weight on a ball in $\RR^d$,
        \textit{Trans. Amer. Math. Soc.}, \textbf{351} (1999),  2439-2458. 
        
\bibitem{X06}
        Yuan Xu,
        A family of Sobolev orthogonal polynomials on the unit ball,
        \textit{J. Approx. Theory}, \textbf{138} (2006), 232-241. 

\end{thebibliography}
\end{document}